\documentclass[english,12pt,twoside]{article}
\usepackage{amssymb,amsbsy,amsmath,amsfonts,amssymb,amscd}
\usepackage{stmaryrd,floatflt}
\usepackage{latexsym}
\usepackage{euscript}
\usepackage{exscale}
\usepackage{epsfig}
\usepackage{amsthm}
\usepackage{babel}

\newtheorem{theo}{Theorem}[section]

\newtheorem{defi}{Definition}[section]
\newtheorem{lem}{Lemma}[subsection]
\newtheorem{prop}{Proposition}[section]
\newtheorem{rem}{Remark}[section]

\newcommand{\K}{\mathfrak K}
\newcommand{\Sp}{\mathrm Sp}

\newcommand{\Hom}{\mathrm{Hom}}

\newcommand{\hatI}{\widehat{\mathfrak I}}
\newcommand{\I}{\mathfrak{I}}

\newcommand{\x}{\underline{x}}
\newcommand{\g}{\underline{g}}
\renewcommand{\k}{\underline{k}}

\newcommand{\h}{\mathfrak h}

\newcommand{\ga}{\gamma}
\newcommand{\Ga}{\Gamma}

\newcommand{\De}{\Delta}

\newcommand{\ph}{\varphi}
\newcommand{\si}{\sigma}

\newcommand{\la}{\lambda}
\newcommand{\La}{\Lambda}

\newcommand{\A}{\mathbb{A}}
\renewcommand{\O}{\mathbb{O}}

\title{Automorphic representations and harmonic cochains for $GL_{n+1}$}
\author{Y. A\"it Amrane}

\begin{document}

\maketitle

\selectlanguage{english} \vspace*{2cm} \hspace{-0.6cm}{\bf
Abstract.} Let $K$ be a global field of positive characteristic. Let
$\infty$ be a fixed place of $K$. This paper gives an explicit
isomorphism between the space of  automorphic forms (resp. cusp
forms) for $GL_{n+1}(K)$ that transform like the special
representations and certain spaces of harmonic cochains (resp. those
with a finite support) defined on the Bruhat-Tits
building of $GL_{n+1}(K_{\infty})$. \bigskip\bigskip

\vspace*{2cm}

\hspace{-0.6cm}{\bf\Large Introduction}\medskip\\
Let $K$ be a global field of characteristic $p>0$, that is a
function field of a smooth geometrically irreducible projective
curve $\cal C$ over a finite field of characteristic $p$. We
identify the places of $K$ with the closed points $|\cal C|$ of
$\cal C$. Let $\infty\in |{\cal C}|$ be a fixed place of $K$ and
$K_\infty$ be the completion of $K$ at this place. Denote by $\A$
and $\A_f$ respectively the ring of ad\`eles and that of finite
ad\`eles of $K$.

Denote by $G$ the reductive group scheme $GL_{n+1}$. Let $M$ be a
commutative ring (with unit) and $L$ be an integral $M$-algebra of
characteristic zero. $M$ and $L$ are assumed to be endowed with the
trivial action of $G(K_{\infty})$. For an open compact subgroup
$\K_f$ of $G(\A_f)$ denote by $\mathfrak{Aut}^{\K_f}(L)$ the space
of $L$-valued automorphic forms that are invariant under ${\mathfrak
K}_f$. These are functions defined on
$$
G(K)\backslash G(\A) /\K_f\times \K_\infty Z_G(K_{_\infty})
$$
for some open compact subgroup $\K_\infty$ of $G(K_\infty)$. Denote
by $\mathfrak{Aut}_{o}^{{\K}_f}(L)$ the subspace of such functions
that moreover are cuspidal.

Let $X$ be a set of representatives of the double cosets
$G(K)\backslash G(\A_f)/\K_f$. It is a finite set. For any $\x\in X$
the intersection $\Ga_{\x}:=G(K)\cap \x \K_f \x^{-1}$ in $G(\A_f)$
is an arithmetic subgroup of $G(K)$. For an integer $k$, $0\leq
k\leq n$, denote by ${\frak Harm}^k(M,L)^{\Ga_{\x}}$ and
respectively ${\frak Harm}_{!}^k(M,L)^{\Ga_{\x}}$ the space of
harmonic cochains on the Bruhat-Tits building of $G(K_\infty)$
invariant under the action of $\Ga_{\x}$  and respectively of those
with finite supports modulo $\Ga_{\x}$.

In this paper, using a result we established in \cite{Yacine} and
that gives in any degree $k$, $0\leq k\leq n$, an explicit
isomorphism between the space of $k$-harmonic cochains and the dual
of the $k$-special representation ${\mbox{Sp}}^{k}(M)$ of
$G(K_{\infty})$, also following the ideas of Drinfeld that did the
work for $G=GL_2$, see \cite{Drinfeld} or also \cite{vdP-R}, we
prove, Th. \ref{lien-frmes-autom-cocycles-harm}, that for every $k$,
$0\leq k \leq n$, we have $M$-isomorphisms :
$$
\bigoplus_{\x \in X} {{\frak Harm}^{k}(M,L)}^{\Ga_{\x}} \cong
{\mbox{Hom}}_{M[G(K_{\infty})]}({\mbox{Sp}}^{k}(M),\; {\frak
Aut}^{\K_{f}}(L))
$$
and
$$
\bigoplus_{\x \in X} {{\frak Harm}_{!}^{k}(M,L)}^{\Ga_{\x}} \cong
{\mbox{Hom}}_{M[G(K_{\infty})]}({\mbox{Sp}}^{k}(M),\; {\frak
Aut}_{\circ}^{\K_{f}}(L)).
$$
These isomorphisms are given explicitly.

We organise this paper as follows. In the first section we give a
brief description of the Bruhat-Tits building associated to
$G(K_{\infty})$. In the second section, we recall the notions of
harmonic cochains and special representations and give the link
between the two. In the third section, we recall the notions of
automorphic forms and cusp forms, we look at them as functions
defined on the ad\`ele groups and we give another interpretation of
them as functions defined on the component at infinity, then we give
other properties in particular when looking at them through the
special representations. Finally, in section 4, we
 give an explicit link between the automorphic and cusp forms that
transform like special representations and the harmonic cochains,
cf. the isomorphisms above.

\section{The Bruhat-Tits building of $GL_{n+1}(K_{\infty})$}

For general properties of buildings, see \cite{Brown} and \cite{Garrett}. An introduction
to the Bruhat-Tits building of $G(K_{\infty})$ with pointed cells is given  in \cite{deShalit}. \medskip \\
{\bf The Bruhat-Tits building (pointed cells).} Let $V$ be the
standard vector space $K_{\infty}^{n+1}$. A lattice in $V$ is a free
$O_{\infty}$-submodule $\La$ of $V$ of rank $n+1$. The Bruhat-Tits
building of $G(K_{\infty})$ may be described as a simplicial complex
$\I$ whose vertices are the dilation classes of lattices. More
precisely, two lattices $\La$ and $\La'$ are in the same class if
$\La'=\la\La$ for some $\la \in K_{\infty}^{^*}$. The class of $\La$
is a vertex $v$ and is denoted $v=[\La]$. For $k$, $0\leq k\leq n$,
a $k$-cell $\si$ in $\I$ is a set of $k+1$ vertices
$\{[\La_{0}],[\La_{1}],\ldots ,[\La_{k}]\}$ such that :
\begin{equation}\label{simplexe}
\cdots \supsetneq \La_{0}\supsetneq \La_{1} \supsetneq \cdots
\supsetneq \La_{k}\supsetneq \pi_{\infty}\La_{0}\supsetneq \cdots
\end{equation}
Notice that there is an obvious cyclic ordering (mod. $(k+1)$) on
the vertices of $\si$.

A pointed $k$-cell of $\I$ is a pair $(\si,v)$ consisting of a
$k$-cell $\si$ together with a distinguished vertex $v$ of $\si$.
Notice, therefore, that in the case of a pointed cell $(\si,v)$
there is a precise ordering on the vertices. If $v=[\La_{0}]$ we
write :
\begin{equation}\label{simplexepointe}
(\si,v) =(\Lambda_{0} \supsetneq \Lambda_{1} \supsetneq \cdots
\supsetneq \Lambda_{k} \supsetneq \pi_{\infty}\Lambda_{0}).
\end{equation}

For each $k$, $0\leq k\leq n$, let $\hatI^{k}$ be the set of pointed $k$-cells of $\I$. \medskip\\
{\bf The action of $G(K_{\infty})$.} For a fixed basis of the vector
space $V$, the action of $G(K_{\infty})$ on $V$ is given by the
matrix product $ug^{-1}$ where $u\in V$ is considered as a line
matrix with respect to the basis of $V$. This action induces an
action of $G(K_{\infty})$ on the vertex set of the building $\I$ by
$$
g.v=[\La g^{-1}].
$$
Thus, $G(K_{\infty})$ acts on the cells by acting on their vertices. \medskip\\
{\bf The type of a pointed cell.} (Cf. \cite[\S\,1.1]{deShalit}.)
Denote by $\kappa_{\infty}=O_\infty/\pi_\infty O_\infty$ the residue
field of $K_\infty$. Let $\si=(\La_{0} \supsetneq \La_{1} \supsetneq
\cdots \supsetneq \La_{k} \supsetneq \pi_{\infty}\La_{0})\in
{\hatI}^{k}$ be a pointed $k$-cell. The type of  $\si$ is defined as
follows :
$$
t(\si)=(d_{1},\ldots ,d_{k+1})
$$
where $d_{i}=\textrm{dim}_{\kappa_{\infty}}\, \La_{i-1}/\La_{i}$ for
each $i=1,\ldots ,k+1$ (here, we suppose
$\La_{k+1}=\pi_{\infty}\La_{0}$).

Since clearly the action of $G(K_{\infty})$ preserves the dimension
of the $\kappa_{\infty}$-vector spaces $\La_{i-1}/\La_{i}$,
 it preserves the type of the pointed $k$-cells as well. \medskip\\
{\bf The standard cells.} Let $\{u_{1},\ldots ,u_{n+1}\}$ be the
standard basis of $V=K_{\infty}^{n+1}$. Consider,  for each
$i=0\ldots n$, the vertex $v^{o}_{i}=[\La_{i}^{^o}]$ represented by
the lattice :
$$
\La^{^o}_{i}=\pi_{\infty} O_{\infty}u_{1}\oplus \cdots \oplus
\pi_{\infty} O_{\infty}u_{i}\oplus O_{\infty}u_{i+1} \oplus \cdots
\oplus O_{\infty}u_{n+1}.
$$
Since the $\La_{i}^{^o},\; 0\leq i\leq n$, satisfy (\ref{simplexe}),
we have an $n$-cell $\si_{\emptyset}=\{v^{o}_{0},v^{o}_{1},\ldots
,v^{o}_{n}\}$ called the fundamental chamber of $\I$.

Now, once and for all, fix $\De=\{1,\ldots ,n\}$. For each
$I\subseteq \De$ such that $\De -I=\{i_{1}< \cdots <i_{k}\}$, we
have a $k$-cell
\begin{equation}\label{simplexe-standard-I}
\si_{I}=\{v^{o}_{0},v^{o}_{i_{1}},\ldots ,v^{o}_{i_{k}}\}.
\end{equation}
The $\si_{I}$, $I\subseteq \De$, are called the standard cells of
the Bruhat-Tits building $\I$. These cells are the faces of the
fundamental chamber $\si_{\emptyset}$ having $v^{o}_{0}$ as vertex,
called the fundamental vertex of $\I$.

We denote by $T$ the standard maximal torus of $G(K_{\infty})$ of
diagonal matrices and by $N$ its normalizer in $G(K_{\infty})$.
Since the Weyl group $W=N/T$ of $G(K_{\infty})$ with respect to $T$
is isomorphic to the permutation group ${\cal S}_{n+1}$, $W$ is
generated by the set $S=\{s_{i},\, i\in \De\}$ of the reflexions
$s_{i}$ which correspond to the transpositions $(i,i+1)\in {\cal
S}_{n+1}$. We have the following lemma :
\begin{lem}\label{wij}
Let $y_{i}$, $0\leq i\leq n$, be the diagonal matrix
$y_{i}=\rm{diag}(\overbrace{1,\ldots
,1}^{i\hbox{\scriptsize~times}},\pi_{\infty},\ldots ,\pi_{\infty})$
and let $w_{i}=(s_{i}s_{i+1}\cdots s_{n})(s_{i-1}s_{i}\cdots
s_{n-1})\cdots (s_{1}s_{2}\cdots s_{n-i+1})\in W$. We have :
$$
(\si_{\emptyset},v_{i}^{o})=y_{i}w_{i}(\si_{\emptyset},v_{0}^{o}).
$$
If $(\si,v_{i_{j}}^{o})=(v_{i_{j}}^{o},\ldots
,v_{i_{k}}^{o},v_{i_{0}}^{o},v_{i_{1}}^{o},\ldots ,v_{i_{j-1}}^{o})$
is a face of the pointed chamber $(\si_{\emptyset},v_{i_{j}}^{o})$,
where $0\leq i_{0}< i_{1}< \cdots < i_{k}\leq n$ and $0\leq j\leq
k$, then
$$
(\si,v_{i_{j}}^{o})=y_{i_{j}}w_{i_{j}}(\si_{{\widehat
I}_{i_{j}}},v_{0}^{o})
$$
where $\De-{\widehat I}_{i_{j}}=\{i_{j+1}-i_{j}<\cdots
<i_{k}-i_{j}<n+1+i_{0}-i_{j}<\cdots <n+1+i_{j-1}-i_{j}\}$.
\end{lem}
\proof The vertices of the fundamental chamber are
$v^{o}_{l}=[\Lambda_{l}^{^o}]$. We can easily check that the
representants $\Lambda^{^o}_{l}$ of these vertices satisfy :
$$
\Lambda^{^o}_{l}y_{i}w_{i}= \left\{ \begin{array}{lll}
\Lambda^{^o}_{n+1+l-i}   &  \textrm{if} & 0\leq l\leq i-1\\
\Lambda^{^o}_{l-i}\pi_{\infty} & \textrm{if} & i\leq l\leq n.
\end{array}
\right.
$$
Therefore, by taking into account the way in which $G(K_{\infty})$
acts on the vertices of $\I$, it follows that :
$$
w_{i}^{-1}y_{i}^{-1}v^{o}_{l}= \left\{ \begin{array}{lll}
v^{o}_{n+1+l-i}   &  \textrm{if} & 0\leq l\leq i-1\\
v^{o}_{l-i} & \textrm{if} & i\leq l\leq n,
\end{array}
\right.
$$
hence
$w_{i}^{-1}y_{i}^{-1}(\si_{\emptyset},v_{i}^{o})=(\si_{\emptyset},v_{0}^{o})$
and, if $(\si,v_{i_{j}}^{o})$ and $\widehat{I}_{i_{j}}$ are as in
the lemma, we have
$w_{i_{j}}^{-1}y_{i_{j}}^{-1}(\si,v^{o}_{i_{j}})=(\si_{\widehat{I}_{i_{j}}},v^{o}_{0})$.
\qed

Since the action of $G(K_{\infty})$ is transitive on the chambers of
$\I$, the lemma above shows that $G(K_{\infty})$ acts transitively
on the pointed $k$-cells of a given type. So, if we denote by
$t_{I}$ the type of the pointed standard $k$-cell
$(\si_{I},v_{0}^{o})$, by $\hatI^{k,t_{I}}$ the set of all pointed
$k$-cells of type $t_{I}$, and by $B_{I}$ the pointwise stabilizer
in $G(K_{\infty})$ of the standard cell $\si_{I}$, we have a
one-to-one correspondence
\begin{equation} \label{pointe}
G(K_\infty)/B_I \xrightarrow{\sim} \hatI^{k,t_{I}}.
\end{equation}\label{transitive}

\section{Harmonic cochains and special representations}

Through all this section, we fix a commutative ring $M$ and an
$M$-module $L$. Assume that $G(K_{\infty})$ acts trivially on $M$
and that $L$ is endowed with an $M$-linear $G(K_{\infty})$-action.

\subsection{Harmonic cochains}

From now on, we sometimes denote by $\si$ a pointed cell $(\si,v)$
when it is clear that it is pointed and which vertex is
distinguished.

Let us recall the definition of harmonic cochains given by E. de
Shalit (\cite[def. 3.1]{deShalit}).

\begin{defi}\label{cocycles-harmoniques}
Let $k$ be an integer such that $0\leq k\leq n$. A $k$-harmonic cochain with
values in the $M$-module $L$ is a homomorphism $\h\in {\Hom}_{M}(M[\hatI^{k}],\,L)$ which
satisfies the following conditions :\\
{\bf (HC1)} If $\si=(v_{0},v_{1},\ldots,v_{k}) \in \hatI^{k}$ is a
$k$-pointed cell and if $\si'=(v_{1},\ldots , v_{k},v_{0})$ is the
same cell but pointed at $v_{1}$, then
$$
\h(\si)=(-1)^{k}\h(\si').
$$
{\bf (HC2)} Fix a pointed $(k-1)$-cell $\eta\in \hatI^{k-1}$, fix a
type $t$ of pointed $k$-cells, and consider the set ${\cal B}(\eta,
t)=\{\si \in \hatI^{k};\, \eta<\si\;{\textrm and}\;t(\si)=t\}$. Then
$$
\sum_{\si\in {\cal B}(\eta,t)}\h(\si)=0.
$$
{\bf (HC3)} Let $k\geq 1$. Fix $\si=(\La_{0}\supsetneq
\La_{1}\supsetneq \cdots \supsetneq \La_{k}\supsetneq
\pi_{\infty}\La_{0})\in \hatI^{k}$ and fix an index $0\leq j\leq k$.
Let ${\cal C}(\si,j)$ be the collection of all $\si'=
(\La'_{0}\supsetneq \La'_{1}\supsetneq \cdots \supsetneq
\La'_{k}\supsetneq \pi_{\infty}\La'_{0})\in \hatI^{k}$ for which
$\La'_{i}=\La_{i}$ if $i\neq j$, $\La_{j}\supseteq
\La'_{j}\supsetneq\La_{j+1}$ and
$\dim_{\kappa}\La'_{j}/\La_{j+1}=1$. Then
$$
\h(\si)=\sum_{\si'\in {\cal C}(\si,j)}\h(\si').
$$
{\bf (HC4)} Let $\si =(v_{0},v_{1},\ldots ,v_{k+1})\in \hatI^{k+1}$.
Let $\si_{j}=(v_{0},\ldots ,{\hat v}_{j},\ldots ,v_{k+1})\in
\hatI^{k}$. Then
$$
\sum_{j=0}^{k}(-1)^{j}\h(\si_{j})=0.
$$
\end{defi}

For any $k$, $0\leq k\leq n$, we denote by ${\mathfrak
Harm}^{k}(M,L)$ the space of $k$-harmonic cochains with values in
the $M$-module $L$. \medskip\\
{\bf The action of $G(K_{\infty})$.} The action of $G(K_{\infty})$
on ${\frak Harm}^{k}(M,L)$ is induced from its natural action on
${\Hom}_{M}(M[\hatI^{k}],\,L)$, namely
$$
(g.{\h})(\si)=g.\h(g^{-1}\si)
$$
for any $\h\in {\mathfrak Harm}^{k}(M,L)$, any $g\in G(K_{\infty})$,
and any $\si \in \hatI^{k}$.

\subsection{Special representations}

Let $P$ be the upper triangular Borel subgroup of $G(K_{\infty})$. A
parabolic subgroup of $G(K_{\infty})$ is a closed subgroup which
contains a Borel subgroup. The subgroups which contain $P$ are said
to be special; these subgroups are completely determined by the
subsets $I$ of $\De$. Indeed, if for each $I\subseteq \De$, we let
$W_{I}$ be the subgroup of $W$ generated by the $s_{i},\, i\in I$,
it has been shown that the subset
$$
P_{I}=P W_{I}P \qquad
(:=\widetilde{P}N_{I}\widetilde{P}\quad\textrm{where}\quad
N_{I}\subseteq N \;\; \textrm{is such that} \;\; N_{I}/T=W_{I})
$$
is a subgroup of $G(K_{\infty})$ containing $P$, and that every
subgroup of $G(K_{\infty})$ containing $P$
is a certain $P_I$ for $I\subseteq \De$. Note that $P=P_{\emptyset}$.\medskip\\
The group $G(K_{\infty})$ is a locally compact topological group,
its topology being induced from that of the non-archimedean field
$K_{\infty}$. We know, that for any $I\subseteq \De$, the
homogeneous space $G(K_{\infty})/P_{I}$ is compact with respect to
the quotient topology.\medskip\\
Denote by $C^{^\infty}(G(K_{\infty})/P_I,M)$ the set of locally
constant functions on  $G(K_{\infty})/P_{I}$ with values in $M$.
\medskip\\
{\bf The action of $G(K_{\infty})$.} The action on
$C^{^\infty}(G(K_{\infty})/P_{I},M)$ is induced by the action by
left translations on
$G(K_{\infty})/P_{I}$.\medskip\\
For any $I_{1}\subseteq I_{2}\subseteq \De$,  we have natural
commutative diagrams of $M[G(K_{\infty})]$-monomorphisms
$$
\begin{array}{lll}
{} & \!\!\!\!\!\!\!\! \!\!\!\!\!\!\!\!C^{^\infty}(G(K_{\infty})/P,\,M)\!\!\!\!\!\!\!\! & {}\\
\qquad\qquad\qquad\;\nearrow & {} &  \!\!\nwarrow \\
C^{^\infty}(G(K_{\infty})/P_{I_{2}},M) & \quad\rightarrow &
\!\!\!\!\!\!\!\!\!C^{^\infty}(G(K_{\infty})/P_{I_{1}},M),
\end{array}
$$
and the special representations of $G(K_{\infty})$ are defined as
follows :
\begin{defi}Let $k$ be an integer with $0\leq k\leq n$ and let $J_{k}$ be the subset
$\llbracket 1,n-k\rrbracket$ of $\De$. A $k$-special representation
of $G(K_{\infty})$ is the $M[G(K_{\infty})]$-module :
$$
\Sp^{k}(M)=\frac{C^{^\infty}(G(K_{\infty})/P_{J_{k}},M)}
{\sum_{j=n-k+1}^{n}C^{^\infty}(G(K_{\infty})/P_{J_{k}\cup\{j\}},M)}.
$$
\end{defi}
In case $k=n$, this is the ordinary Steinberg representation.

One thing important to know about the special  representations is
that they are cyclic. Indeed : \\

For each $I\subseteq \De$, we denote by $B_{I}^{^\circ}$ the open
compact subgroup of $G(O_{\infty})$ which is the inverse image of
the standard parabolic subgroup $P_{I}(\kappa_\infty)$ of
$G(\kappa_\infty)$ by the map ``reduction mod. $\pi_{\infty}$'' :
$G(O_{\infty})\rightarrow G(\kappa_\infty)$. The parahoric subgroups
of $G(K_{\infty})$ are  the conjugates in $G(K_{\infty})$ of the
$B_{I}^{^\circ}$, $I\subseteq \De$. Note that we have
\begin{equation}\label{BI}
B_{I}=B_{I}^{^\circ}K_{\infty}^{*},
\end{equation}\label{parahoric-subgroups}
and that then the $B_I$ is compact open modulo the center of
$G(K_\infty)$.

Let $I\subseteq \De$. For any subset $H$ of $G(K_{\infty})$, we
denote by $\chi_{HP_{I}}\in C^{^\infty}(G(K_{\infty})/P_{I},M)$ the
characteristic function of $HP_{I}/P_{I}$.

\begin{prop}\label{chibp} (P. Schneider and U. Stuhler)
The $M[G(K_{\infty})]$-module $C^{^\infty}(G(K_{\infty})/P_{I},M)$
is generated by the characteristic function $\chi_{B_{I}P_{I}}$.
\end{prop}
\proof See \cite[\S 4, prop. 8' and cor. 9']{Schneider}. \qed

\subsection{Harmonic cochains and special representations}
Let $k$ be an integer, $0\leq k\leq n$. In \cite{Yacine} we proved
that the space of harmonic cochains of degree $k$ is isomorphic to
the dual of the $k$-special representation. For this purpose,
inspired by the definition of the harmonic cochains and using
parahoric subgroups, we defined a certain type of sets. Here we
recall these essential points.

Let $I\subseteq \De$ such that $\De-I=\{i_1 <i_2<\cdots <i_k\}$. It
is clear that for any $m=1, \ldots,k$, we have $i_m\leq n-k+m$. Set
\begin{equation}\label{CI}
C_{I}^{^\circ}=B^{^\circ}_{\llbracket i_k+1,n \rrbracket}
B^{^\circ}_{\llbracket i_1+1,n-k+1 \rrbracket}B^{^\circ}_{\llbracket
1,n-k \rrbracket}\quad \textrm{and} \quad C_{I}=B_{\llbracket
i_k+1,n \rrbracket} B_{\llbracket i_1+1,n-k+1
\rrbracket}B_{\llbracket 1,n-k \rrbracket}
\end{equation}\label{defCI}
The set $C_I^{^\circ}$ is compact open in $G(K_{\infty})$ and we
have $C_I = C_I^{^\circ} K^{*}_{\infty}$, see
\S\ref{parahoric-subgroups}. Hence, the set
$C_IP_{J_k}/P_{J_k}=C^{^\circ}_IP_{J_k}/P_{J_k}$ is compact open in
$G(K_{\infty})/P_{J_k}$.

\begin{theo}\label{maintheorem} For any $k$, $0\leq k\leq n$, there is an $M[G(K_{\infty})]$-isomorphism
$$
{\mathfrak Harm}^{k}(M,L) \cong {\Hom}_{M}(\Sp^{k}(M),\,L)
$$
\end{theo}
\proof The proof here involves hard computations of combinatorial
nature. The interested reader could find it in \cite{Yacine}. For
later use we just need to describe the isomorphism. The isomorphism
is given by the map :
$$
\Phi : {\mathfrak Harm}^{k}(M,L) \longrightarrow
{\Hom}_{M}(\Sp^{k}(M),\,L)
$$
which to a harmonic cochain $\h$ associates the $M$-linear map
 $\ph_{\h}$ defined by
\begin{equation}\label{phi-h}
\ph_{\h}(g\chi_{B_{J_k}P_{J_k}})=\h(g(\si_{J_k},v_{0}^{o})),
\end{equation}
for any $g\in G(K_{\infty})$. The inverse map
$$
H : {\Hom}_{M}(\Sp^{k}(M),\,L) \longrightarrow  {\mathfrak
Harm}^{k}(M,L)
$$
sends $\ph$ to $\h_{\ph}$ given by
\begin{equation}\label{h-phi}
\h_{\ph}(g(\si_{I},v_{0}^{o}))=\ph(g\chi_{C_{I}P_{J_k}}),
\end{equation}
for any $g\in G(K_{\infty})$. \qed

\section{Automorphic forms}

\subsection{Automorphic forms}

Let $\A=\prod_{\nu \in |{\cal C}|}'K_{\nu}$ be the ring of ad\`eles
of $K$, i.e. the restricted product of the family  $(K_{\nu})_{\nu
\in |X|}$ with respect to the family of the compact open subrings
$(O_{\nu})_{\nu \in |{\cal C}|}$ :
$$
\mathbb{A}=\{(a_{\nu})\in \prod_{\nu \in |{\cal C}|}K_{\nu} |\;
a_{\nu}\in O_{\nu}\;\textrm{for almost all}\;\nu\in |{\cal C}|\}
$$

\par We denote by $\O=\prod_{\nu\in |{\cal C}|}O_{\nu}$, it is an open compact
subset of $\mathbb A$. We can write $\A=\A_{f}\times K_{\infty}$ et
$\O=\O_f\times O_{\infty}$, where we have set $\A_f=\prod_{\nu \in
|{\cal C}|-\{\infty\}}'K_{\nu}$ the ring of finite ad\`eles of $K$
and $\O_{f}=\prod_{\nu\in |{\cal C}|-\{\infty\}}O_{\nu}$.

\par For any $a \in K$, we have $\nu (a)= 0,\; \mbox{for almost all}\; \nu \in |{\cal C}|$, which means that $K$ can
be seen as a subfield of $\mathbb A$ imbedded diagonally.

\par In what follows we work with  $G(K),G(\A)$ et $G(\O)$ \ldots that have
the same properties than $K,\A ,\O \ldots$ recalled above. For
example, $G(K)$ is embedded diagonally in $G(\A)$, $G(\O)$ is a
compact open subgroup of $G(\A)$ \ldots

Let $M$ be a commutative ring with a unit and  $L$ be an integral
$M$-algebra of characteristic zero. Assume these are endowed with
the trivial action of $G(\A)$.

\begin{defi}\label{formes-automorphes}(G. Harder, \cite{Harder}) An
automorphic form with values in $L$, with respect to an open compact
subgroup $\K$ of $G(\O)$, is a function $f:G(\A)\rightarrow
L$, such that $f(\ga \g\k) =f(\g)$ for  $\ga \in G(K)$, $\g \in
G(\A)$ and $\k\in \K Z_{G}(K_{\infty})$.\\
An automorphic form is a cusp form if moreover :
$$
\int_{U_{I}(K)\backslash
U_{I}(\A)}f(\underline{u}\g)\textrm{d}\underline{u}_{I} =0
$$
for every $I\subseteq \Delta$, ($\textrm{d}\underline{u}_{I}$ a Haar
measure which is normalized with respect to the compact
$U_{I}(K)\backslash U_{I}(\A)$).
\end{defi}
In fact this integral is a finite sum, so it has a sense for any
choice of $L$ of characteristic zero.

Denote by ${\mathfrak Aut}(\K,L)$ (resp. ${\mathfrak
Aut}_{0}(\K,L)$) the set of L-valued automorphic forms with respect
to an open compact subgroup $\K$ of $G(\O)$ (resp. the set of those that are cusp forms).\\

\par From now on, once and for all we fix an open compact subgroup $\K_f$ of $G(\O_{f})$,
 $X=X_{\K_{f}}$ a set of representatives of the double cosets $G(K)\backslash
G(\A_{f})/\K_{f}$. For every $\x\in X$, we let $\Ga_{\x}$ the
arithmetic group ${\x}\K_{f}{\x}^{-1}\cap G(K)$. Put
$$
{\mathfrak Aut}^{\K_f}(L)=\bigcup_{\K_{\infty}}{\mathfrak Aut}(\K_f
\times \K_{\infty},\, L) \;\;\;\mbox{and}\;\;\; {\mathfrak
Aut}_0^{\K_f}(L)=\bigcup_{\K_{\infty}}{\mathfrak Aut}_0(\K_f\times
\K_{\infty},\, L)
$$
where  $\K_{\infty}$ runs through the open compact subgroups of $G(O_{\infty})$.\medskip\\
$\bullet$ {\bf The action $G(K_{\infty})$ :} $G(K_{\infty})$ acts on
${\mathfrak Aut}^{\K_{f}}(L)$ by the formula
\begin{equation}
(g.f)(\g)=f(\g g).
\end{equation}
Indeed, it is not difficult to see that if $f$ is in ${\mathfrak
Aut}(\K_{f}\times \K_{\infty},\, L)$ then $g.f$ belongs to
${\mathfrak Aut}(\K_{f}\times g\K_{\infty}g^{-1}\cap
G(O_{\infty}),\, L) \subseteq
{\mathfrak Aut}^{\K_{f}}(L)$. \\

We end this paragraph by stating the following important result
about cusp forms :

\begin{theo}(G. Harder)\label{harder}
Let $\mathfrak K$ be an open compact subgroup of $G(\A)$. Then there
 exists an open subset ${\cal U}_\K$ of $G(\A)$ such that
 $G(K){\cal U}_\K{\frak K}Z_{G}(K_{\infty})={\cal U}_\K$, the
quotient $G(K)\backslash {\cal U}_\K/\K Z_{G}(K_{\infty})$ is
finite, and we have :
$$
\rm{supp}(f) \subseteq {\cal U}_\K, \qquad {\rm\it pour\; tout}\; f
\in {\mathfrak Aut}_{0}(\K,L).
$$
\end{theo}
\proof See \cite[cor. 1.2.3]{Harder}. \qed

\subsection{Automorphic forms as functions on $G(K_\infty)$}

\begin{defi}
A subgroup $\Ga$ of $G(K)$ is said to be arithmetic if it is
commensurable with $G(A)$, i.e. if $\Ga \cap G(A)$ is a subgroup of
finite index in both $\Ga$ and $G(A)$.
\end{defi}

\par For any subgroup ${\mathfrak H}_{f}$ of $G(\A_f)$, choose a set $X_{{\mathfrak H}_{f}}\subseteq G(\A_f)$
of representatives of the double cosets $G(K)\backslash
G(\A_{f})/{\mathfrak H}_{f}$. For every $\x \in X_{{\mathfrak
H}_{f}}$, let
$$
\Ga_{\x}=G(K)\cap {\x}{\mathfrak H}_{f}{\x}^{-1}
$$
be the intersection in $G(\A_{f})$. It is a discrete subgroup of
$G(K_{\infty})$.

The following is a well known result, see \cite{Harder}.

\begin{prop}
If a subgroup ${\mathfrak H}_{f}$ of $G(\A_{f})$ contains an open
compact subgroup of $G(\A_{f})$, then $X_{{\mathfrak H}_{f}}$ is a
finite set. If, moreover, ${\mathfrak H}_{f}$ is open compact in
$G(\A_{f})$, the $\Ga_{\x}$ ($\x\in X_{{\mathfrak H}_{f}}$) are
arithmetic.
\end{prop}

We have :

\begin{lem}\label{adinfini}
Let $H_{\infty}$ be a subgroup of $G(K_{\infty})$ and ${\frak
H}_{f}$ be a subgroup of $G(\A_{f})$. For any  $\x\in X_{{\frak
H}_{f}}$, let $\Ga_{\x}=\x{\mathfrak H}_{f}{\x}^{-1}\cap G(K)$. We
have a bijective correspondence  :
$$
G(K)\backslash G(\A) /({\mathfrak H}_{f}\times H_{\infty})
\xrightarrow{\sim} \coprod_{\x\in X_{{\mathfrak H}_{f}}}
\Ga_{\x}\backslash G(K_{\infty})/H_{\infty}
$$
given by the map which sends $G(K)\g({\mathfrak H}_{f}\times
H_{\infty})$ to the double class
$\Ga_{\x}\,\tau^{-1}g_{\infty}H_{\infty}$, for any $\g=(\tau \x
\underline{h}_{f},g_{\infty}) \in G(\A)$.
\end{lem}
\proof $\;$\smallskip\\
The map is well defined. Indeed, let $\g=(\g_{f},g_{\infty})\in
G(\A)$ and let $\tau,\tau' \in G(K)$ and
$\underline{h}_{f},\underline{h}'_{f}\in {\mathfrak H}_{f}$ be such
that $\g_{f}=\tau\x\underline{h}_{f}=\tau'\x\underline{h}'_{f}$.
This implies that
$\tau^{-1}\tau'=\x\underline{h}_{f}{\underline{h}'_{f}}^{-1}\x^{-1}
\in \x{\mathfrak H}_{f}\x^{-1}\cap G(K)=\Ga_{\x}$, then :
$$
\tau^{-1}g_{\infty}=(\tau^{-1}\tau'){\tau'}^{-1}g_{\infty}\in
\Ga_{\x}{\tau'}^{-1}g_{\infty}.
$$
For the infinite component, The right invariance of the map
$H_{\infty}$ is clear. \medskip\\
The map is surjective. For any $\x \in X_{{\mathfrak H}_{f}}$ and
any $g_{\infty}\in G(K_{\infty})$, the double coset
$\Ga_{\x}g_{\infty}H_{\infty}$ is the image of
$G(K)(\x,g_{\infty})({\mathfrak H}_{f}\times H_{\infty})$.\medskip\\
The map is injective. Let $\g=(\tau\x\underline{h}_{f},g_{\infty})$
and $\g'=(\tau'\x'\underline{h}_{f}',g'_{\infty})$ be such that
$\Ga_{\x}\tau^{-1}g_{\infty}H_{\infty}=\Ga_{\x'}\tau'^{-1}g'_{\infty}H_{\infty}$.
The union over $X_{{\mathfrak H}_{f}}$ being a disjoint union, we
must have $\x=\x'$. Thus, if $\ga \in \Ga_{\x}$ and $h_{\infty}\in
H_{\infty}$ are such that $\tau^{-1}g_{\infty}=\ga
\tau'^{-1}g'_{\infty}h_{\infty}$, we have
$$
\g=(\tau\x\underline{h}_{f},g_{\infty})=(\tau\x\underline{h}_{f},\tau\ga\tau'^{-1}g'_{\infty}h_{\infty})
=\tau\ga\tau'^{-1}(\tau'\ga^{-1}\x,g'_{\infty})(\underline{h}_{f},h_{\infty}).
$$
Now, since $\ga^{-1} \in \Ga_{\x}=\x{\mathfrak H}_{f}\x^{-1}\cap
G(K)$, there exists $\underline{h}''_{f}\in {\mathfrak H}_{f}$ such
that $\ga^{-1} =\x\underline{h}''_{f}\x^{-1}$. So
$$
\g=\tau\ga\tau'^{-1}(\tau'\x,,g'_{\infty})(\underline{h}''_{f}\underline{h}_{f},h_{\infty})
=\tau\ga\tau'^{-1}\g'(\underline{h}_{f}'^{-1}\underline{h}''_{f}\underline{h}_{f},h_{\infty}),
$$
which means that $\g$ and $\g'$ are representatives of the same
double coset modulo $G(K)$ on the left and modulo ${\mathfrak
H}_{f}\times H_{\infty}$ on the right. \qed

This lemma \ref{adinfini} gives us a way to interpret the
automorphic forms as functions defined on the quotient sets
$\Ga_{\x}\backslash G(K_{\infty})$ right invariant under some open
compact subgroup $\K_{\infty}$ of $G(O_{\infty})$ and under
$Z_{G}(K_{\infty})$. More precisely, if for each $\x\in X$ one
defines :
$$
{\mathfrak Aut}_{\x}(L) \qquad (\textrm{resp.} \qquad {\frak
Aut}_{\x,\circ}(L))
$$
to be the set of functions $f: G(K_{\infty}) \shortrightarrow L$,
left invariant under $\Ga_{\x}$ and right invariant under
 $\K_{\infty}Z_{G}(K_{\infty})$, where ${\K}_{\infty}$ is an open compact subgroup of $G(O_{\infty})$ which
 depends on $f$ (respectively, the set of such functions that
moreover have a finite support in $\Ga_{\x}\backslash
G(K_{\infty})/{\K}_{\infty}Z_{G}(K_{\infty})$). Then we have the
following :
\begin{prop}\label{frm.aut.infty} We have isomorphism of
$M[G(K_{\infty})]$-modules :
$$
\Xi : \displaystyle \bigoplus_{\x \in X} {\mathfrak Aut}_{\x}(L)
\xrightarrow{\cong}{\mathfrak Aut}^{\K_f}(L)
$$
which maps $(f_{\x})_{\x}$ to $f$ defined by the formula $
f(\g)=f_{\x}(\tau^{-1}g_{\infty})$ for $\g=(\tau
\x\k_f,g_{\infty})\in G(\A)$. Moreover, by this isomorphism when we
restrict to the cusp forms we get an isomorphism of
$M[G(K_{\infty})]$-modules :
$$
\displaystyle \bigoplus_{\x \in X}{\mathfrak
Aut}_{\x,\circ}(L)\cong{\mathfrak Aut}_{\circ}^{\K_f}(L)
$$
\end{prop}
\proof The first isomorphism is a direct consequence of Lemma
\ref{adinfini}. To get the second isomorphism we only need to use
Theorem \ref{harder}.  \qed

We finish this section with the following remark that motivates out
link between automorphic forms and harmonic cochains.

\begin{rem} Let us for while get back to the Bruhat-Tits building. For any
$I\subseteq \Delta$ such that $|\Delta - I|=k$, if in the lemma
above we put $H_{\infty}=B_{I}=B_{I}^{^\circ}Z_{G}(K_{\infty})$,
then, see (\ref{pointe}) in \S\,\ref{transitive}, we get a
one-to-one correspondence :
$$
G(K)\backslash G(\A) /({\mathfrak H}_{f}\times
B_{I}^{^\circ}Z_{G}(K_{\infty})) \xrightarrow{\sim} \coprod_{\x\in
X} \Ga_{\x}\backslash {\widehat{\mathfrak I}}^{k,t_{I}}.
$$
In this section the link between automorphic forms, functions
defined on
 $G(K)\backslash G(\A)$ and right invariant by an open compact subgroup
 of $G(\mathbb{A})$ and by $Z_{G}(K_{\infty})$, and harmonic cochains defined
 on ${\widehat{\mathfrak I}}^{k}$ and which are invariant under the action of the arithmetic groups $\Ga_{\x}$, $\x\in X$.
\end{rem}

\subsection{Automorphic forms through special representations.}

Recall, see Proposition \ref{chibp}, that as an
$M[G(K_{\infty})]$-module the space
$C^{\infty}(G(K_{\infty})/P_{J},\, M)$ is generated by the
characteristic function $\chi_{B_{J}P_{J}}$. So every element
$\varrho\in {\rm Sp}^{k}(M)$ is of the form :
$$
\varrho=\sum_{j=1}^{m}\alpha_{j}\varrho_{\infty}^{j}\chi_{B_{J}P_{J}}
$$
with $\alpha_{j}\in M$ and $\varrho_{\infty}^{j}\in G(K_{\infty})$.
Therefore by the action of $G(K_{\infty})$  this $\varrho$ is
invariant under the open compact subgroup :
\begin{equation}\label{kh}
{\frak
K}_{\infty}^{\varrho}=\bigcap_{j=1}^{m}\varrho^{j}_{\infty}B_{J}^{^\circ}
{\varrho_{\infty}^{j}}^{\!\!\!-\!1}\cap G(O_{\infty})
\end{equation}
of $G(O_{\infty})$.

\begin{prop}
Let $\cal F$ be the set of all functions $f:G({\mathbb
A}_{f})\rightarrow  L$ endowed with the action of $G(K)$ coming from
 left translation on $G({\mathbb A}_{f})$. We have an
 $M$-isomorphism
 $$
 \Psi:\; \mbox{Hom}_{M}({\rm Sp}^{k}(M),\, {\cal
F})^{G(K)} \xrightarrow{\cong} \mbox{Hom}_{M[G(K_{\infty})]}({\rm
Sp}^{k}(M),\,{\frak Aut}^{\K_f}(L)).
$$
which sends $\ph$ to $\psi_{\ph}$ defined as follows. For any
$\varrho \in {\rm Sp}^{k}(M)$, the function
$\psi_{\ph}(\varrho):G(\A)\shortrightarrow L$ is given by the
formula
$$
\psi_{\ph}(\varrho)(\g)=\ph(g_{\infty}. \varrho)(\g_f)
$$
for $\g=(\g_f,g_{\infty})\in G(\A)$.
\end{prop}
\proof $\;$\\
$\Psi$ is well defined. Indeed, let us show first that for any
$\varrho \in {\rm Sp}^{k}(M)$, we have $\psi_{\ph}(\varrho)\in
{\frak Aut}^{\K_f}(L)$. Recall that $G(K)$ is seen as a subgroup of
$G(\A_f)$ and of $G(\A)$ embedded diagonally, so let $\tau \in G(K)$
and write $\tau=(\underline{\tau}_f,\tau_{\infty})$ where
$\underline{\tau}_{f}=\tau \in G({\A}_{f})$ and $\tau_{\infty}=\tau
\in G(K_{\infty})$. Let $\K^{\varrho}_{\infty}$ be the open compact
subgroup of $G(O_{\infty})$ given by (\ref{kh}) and recall that then
$\varrho$ is invariant under
$\K^{\varrho}_{\infty}Z_{G}(K_{\infty})$. Let $\g=(\g_f,
g_{\infty})\in G(\A)$ and $\k=(\k_{f}, k_{\infty})\in \K_f \times
\K^{\varrho}_{\infty}Z_{G}(K_{\infty})$. We have :
$$
\begin{array}{ll}
\psi_{\ph}(\varrho)(\tau \g\k) &
=\psi_{\ph}(\varrho)(\underline{\tau}_{f}{\g}_{f}{\k_{f}},\tau_{\infty}g_{\infty}k_{\infty})\\
& =[\ph(\tau_{\infty}g_{\infty}k_{\infty}.\varrho)](\underline{\tau}_{f}{\g}_{f}{\k_{f}})\\
& =[({\tau}^{-1}.\ph)(g_{\infty}k_{\infty}.\varrho)](\g_f\k_f)\\
& =\ph(g_{\infty}.\varrho)(\g_f)\\
& =\psi_{\ph}(\varrho)(\g).
\end{array}
$$
The equality before the last comes from the invariance of $\ph$
under $G(K)$, of $\varrho$ under
$\K^{\varrho}_{\infty}Z_G(K_{\infty})$ and of $\cal F$ on the right
under $\K_f$. Therefore $\psi_{\ph} \in {\frak Aut}(\K_f\times
\K^{\varrho}_{\infty},L)\subseteq {\frak Aut}^{\K_f}(L)$. Now let us
check that $\psi$ is a homomorphism of $M[G(K_{\infty})]$-modules.
For any $\varrho\in {\rm Sp}^{k}(M)$, $\g=(\g_f,g_{\infty})\in
G(\A)$ and $u_{\infty}\in G(K_\infty)$, we have :
$$
\begin{array}{ll}
\psi_{\ph}(u_{\infty}.\varrho)(\g) & =\ph(g_{\infty}u_{\infty}.\varrho)(\g_f)\\
&=\psi_{\ph}(\varrho)(\g_f,g_{\infty}u_{\infty})\\
&=\psi_{\ph}(\varrho)(\g u_{\infty})\\
&=[u_{\infty}.\psi(\varrho)](\g),
\end{array}
$$
therefore
$\psi_{\ph}(u_{\infty}.\varrho)=u_{\infty}\psi_{\ph}(\varrho)$ for
any $\varrho\in {\rm Sp}^{k}(M)$ and any $u_{\infty}\in
G(K_{\infty})$. Consequently, $\psi_{\ph}$ is a homomorphism of
$M[G(K_{\infty})]$-modules and $\Psi$ is well defined clearly
$M$-linear map.

To prove that $\Psi$ is an isomorphism we give its reciprocal
$$
\begin{array}{cccc} \Psi^{'} : & {\rm
Hom}_{M[G(K_{\infty})]}({\rm Sp}^{k}(M),\,{\frak Aut}^{\K_f}(L)) &
\longrightarrow & {\rm Hom}_{M}({\rm Sp}^{k}(M),\, {\cal
F})^{G(K)}\\ {} & \psi & \longmapsto & \ph_{\psi} : {\rm Sp}^{k}(M)
\mapsto {\cal F}
\end{array}
$$
where for any $\varrho\in {\rm Sp}^{k}(M)$, the function $\ph_{\psi}
: G(\A_{f})\shortrightarrow L$ is given by the formula :
$$
\ph_{\psi}(\varrho)(\g_f)=\psi(\varrho)((\g_f,1_{\infty}))
$$
for $\g_f\in G(\A_f)$.\\
Let us check that indeed $\ph_{\psi}(\varrho)\in {\cal F}$ for
$\varrho \in {\rm Sp}^{k}(M)$. For any $\g_f\in G(\A_f)$ and any
$\k_f\in \K_f$, we have :
\begin{equation}\label{e4}
\ph_{\psi}(\varrho)(\g_f\k_f) =\psi(\varrho)(\g_f\k_f,1_{\infty})
=\psi(\varrho)((\g_f,1_{\infty})(\k_f,1_{\infty})),
\end{equation}
but $\psi(\varrho)\in {\frak Aut}^{\K_f}(L)$, there is then an open
compact subgroup $\K_{\infty}$ of $G(O_{\infty})$ such that
$\psi(\varrho)$ is right invariant under $\K_f\times \K_\infty$. We
have $(\k_f,1_{\infty}) \in \K_f\times \K_\infty$, hence :
$$
\psi(\varrho)((\g_f,1_{\infty})(\k_f,1_{\infty}))=\psi(\varrho)(\g_f,1_{\infty})
=\ph_{\psi}(\varrho)(\g_f);
$$
these equalities with that of (\ref{e4}) give
$\ph_{\psi}(\varrho)(\g_f\k_f)=\ph_{\psi}(\varrho)(\g_f)$ and therefore $\ph_{\psi}(\varrho)\in {\cal F}$.\\
Let us now check that $\ph_{\psi}$ is invariant under the action
$G(K)$. Let $\tau\in G(K)$,
$\tau=(\underline{\tau}_{f},\tau_{\infty})\in G(\A_f)\times
G(K_{\infty})$. By the way $\tau$ acts on $\ph_\psi$  and on the
elements of $\cal F$, for any $\varrho\in {\rm Sp}^{k}(M)$ and any
$\g_f\in G(\A_f)$, we have :
\begin{equation}\label{blabla1}
[(\tau.\ph_psi)(\varrho)](\g_f)=[\underline{\tau}_{f}.\ph_\psi(\tau_{\infty}^{-1}.\varrho)](\g_f)
=\ph_\psi(\tau_{\infty}^{-1}.\varrho)(\underline{\tau}_{f}^{-1}\g_f)
=\psi(\tau_{\infty}^{-1}.\varrho)(\underline{\tau}_{f}^{-1}\g_f,1_{\infty}),
\end{equation}
than, since $\psi$ is a homomorphism of $M[G(K_{\infty})]$-modules,
we have :
\begin{equation}\label{blabla2}
\psi(\tau_{\infty}^{-1}.\varrho)(\underline{\tau}_{f}^{-1}\g_f,1_{\infty})
=[\tau_{\infty}^{-1}.\psi(\varrho)](\underline{\tau}_{f}^{-1}\underline{g}_{f},1_{\infty})
\end{equation}
and by the way $G(K_{\infty})$ acts on $\psi(\varrho)\in {\frak
Aut}^{\K_f}(L)$, we have :
\begin{equation} \label{blabla2.5}
[\tau_{\infty}^{-1}.\psi(\varrho)](\underline{\tau}_{f}^{-1}\underline{g}_{f},1_{\infty})
=\psi(\varrho)((\underline{\tau}_{f}^{-1}\underline{g}_{f},1_{\infty})\tau_{\infty}^{-1})
=\psi(\varrho)(\underline{\tau_{f}}^{-1}\underline{g}_{f},\tau_{\infty}^{-1}),
\end{equation}
and finally, as $\tau$ seen in $G(\A)$ is equal to
$(\tau_{f},\tau_{\infty})$, we have
$(\underline{\tau_{f}}^{-1}\g_f,\tau_{\infty}^{-1})=\tau^{-1}(\g_f,1_{\infty})$,
hence :
\begin{equation}\label{blabla3}
\psi(\varrho)(\underline{\tau_{f}}^{-1}\g_f,\tau_{\infty}^{-1})
=\psi(\varrho)(\tau^{-1}(\g_f,1_{\infty}))
=\psi(\varrho)(\g_f,1_{\infty}) =\ph_\psi(\varrho)(\g_f).
\end{equation}
From the equalities (\ref{blabla1}), (\ref{blabla2}),
(\ref{blabla2.5}) and (\ref{blabla3}) we conclude that we have
$[(\tau.\ph_\psi)(\varrho)](\g_f)=\varphi(\varrho)(\g_f)$, therefore
$\ph_\psi$ is invariant under the action of $\tau \in G(K)$.

The $M$-linear maps $\Psi$ et $\Psi'$ are inverse of each other. On
one hand we have $\Psi \circ \Psi'=I$. Indeed, let $\psi\in {\rm
Hom}_{M[G(K_{\infty})]}({\rm Sp}^{k}(M),\,{\frak Aut}^{\K_f}(L))$,
let $\ph_{\psi}$ be its image by $\Psi'$ and let $\psi_{\ph_{\psi}}$
 be the image of $\ph_{\psi}$ by $\Psi$. Thus, for any $\varrho\in
{\rm Sp}^{k}(M)$ and any $\g=(\g_f,g_{\infty})\in G(\A)$, we have :
\begin{equation}\label{e1}
\psi_{\varphi_{\psi}}(\varrho)(\underline{g})
=\varphi_{\psi}(g_{\infty}.\varrho)(\underline{g}_{f})
=\varphi_{\psi}(g_{\infty}.\varrho)(\underline{g}_{f})
=\psi(g_{\infty}.\varrho)(\underline{g}_{f},1_{\infty}),
\end{equation}
and $\psi$ being a homomorphism of $M[G(K_{\infty})]$-modules, also
by the $G(K_{\infty})$ acts on the automorphic form $\psi(\varrho)$,
we have :
$$
\psi(g_{\infty}.\varrho)(\g_f,1_{\infty})
=(g_{\infty}.\psi(\varrho))(\g_f,1_{\infty})
=\psi(\varrho)((\g_f,1_{\infty})g_{\infty})
=\psi(\varrho)(\g_f,1_{\infty}g_{\infty}) =\psi(\varrho)(\g).
$$
Combining these equalities with that of (\ref{e1}), we deduce that
we have $\psi_{\ph_{\psi}}(\varrho)(\g)=\psi(\varrho)(\g)$ for any
$\varrho\in {\rm Sp}^{k}(M)$ and any $\underline{g}\in G({\mathbb
A})$, hence $\psi_{\varphi_{\psi}}=\psi$.\\
On the other hand we have $\Psi' \circ \Psi =I$. Indeed, let $\ph
\in {\rm Hom}_{M}({\rm Sp}^{k}(M),\, {\cal F})^{G(K)}$, let
$\psi_{\ph}=\Psi(\ph)$ and $\ph_{\psi_{\ph}}=\Psi'(\psi_{\ph})$. As
above, we need to prove that $\ph_{\psi_{\ph}}=\ph$. For any
$\varrho\in {\rm Sp}^{k}(M)$ and any $\g_f\in G(\A_f)$, we have :
$$
\ph_{\psi_{\ph}}(\varrho)(\g_f)
=\psi_{\ph}(\varrho)(\g_f,1_{\infty}) =\ph(1_{\infty}.\varrho)(\g_f)
=\ph(\varrho)(\g_f).
$$
Therefore, $\Psi$ is an $M$-isomorphism. \qed

\begin{prop}\label{frm.aut.repr.spec.} For every integer $k$, $0\leq
k \leq n$. For each $\x\in X$ , we have an isomorphism of
$M$-modules :
$$
\Psi_{\x}:\; \mbox{Hom}_{M}({\rm Sp}^{k}(M),\,
L)^{\Gamma_{\x}}\xrightarrow{\cong}
\mbox{Hom}_{M[G(K_{\infty})]}({\rm Sp}^{k}(M),\,{\frak
Aut}_{\x}(L)).
$$
which to $\ph$ associates $\psi_{\ph}$ defined as follows. For any
$\varrho\in {\rm Sp}^{k}(M)$, the function $\psi_{\ph}(\varrho) :
G(K_{\infty})\rightarrow L$ is given by the formula :
$$
\psi_\ph(\varrho)\,(g_{\infty})=\varphi(g_{\infty}.\varrho).
$$
\end{prop}

\proof $\;$\\
First, we need to prove that $\Psi_{\x}$ is well defined. That is
$\psi_{\ph}(\varrho) \in {\frak Aut}_{\x}(L)$ and $\psi_{\ph}$ is an
$M[G(K_{\infty})]$-homomorphism. Since $\varphi : {\rm Sp}^{k}(M)
\shortrightarrow L$ is an $M$-linear map invariant under $\Ga_{\x}$,
and $\varrho$ is invariant under ${\frak
K}_{\infty}^{\varrho}Z_{G}(K_{\infty})$ where ${\frak K}^{\varrho}$
is the compact open subgroup given by (\ref{kh}) above, then for any
$\ga \in \Ga_{\x}$, $g_{\infty}\in G(K_{\infty})$ and $k_{\infty}\in
{\frak K}^{\varrho}_{\infty}Z_{G}(K_{\infty})$, we have
$$
\psi_{\ph}(\varrho)\,(\gamma g_{\infty}k_{\infty})=\ph (\ga
g_{\infty}k_{\infty}.\varrho)=(\ga^{-1}.\ph)(g_{\infty}.(k_{\infty}.\varrho))
=\ph(g_{\infty}.\varrho)=\psi_{\ph}(\varrho)(g_{\infty}).
$$
It is clear that $\psi_{\ph}$ is $M$-linear, so it remains to check
that it is $G(K_{\infty})$-equivariant. Indeed it is, we have :
$$
\psi_{\ph}(u_{\infty}.\varrho)\,(g_{\infty})
=\ph(g_{\infty}.(u_{\infty}.\varrho))
=\ph((g_{\infty}u_{\infty}).\varrho)
=\psi_{\ph}(\varrho)(g_{\infty}u_{\infty})
=[u_{\infty}.\psi_{\ph}(\varrho)](g_{\infty}).
$$
for any $g_{\infty}\in G(K_{\infty})$, any $u_{\infty} \in
G(K_{\infty})$ and any $\varrho\in {\rm Sp}^{k}(M)$. \\

Now, in order to prove that $\Psi_{\x}$ is an isomorphism, we give
its reciprocal map. To an $M[G(K_{\infty})]$-homomorphism $\psi
:{\rm Sp}^{k}(M) \shortrightarrow {\frak Aut}_{\x}(L)$ we associate
a map $\ph_{\psi} :{\rm Sp}^{k}(M) \shortrightarrow L $ given by the
formula :
$$
\ph_{\psi}(\varrho)=\psi(\varrho)(1_{\infty}),
$$
for any $\varrho \in {\rm Sp}^{k}(M)$ and where $1_{\infty}$ is the
identity element of $G(K_{\infty})$. We need only to check that
$\ph_{\psi}$ is invariant under $\Gamma_{\x}$ to get clearly an
$M$-linear map :
$$
\Psi_{\x}^{'}:\; {\rm Hom}_{M[G(K_{\infty})]}({\rm
Sp}^{k}(M),\,{\frak Aut}_{\x}(L))\longrightarrow {\rm Hom}_{M}({\rm
Sp}^{k}(M),\, L)^{\Gamma_{\x}}.
$$
For any $\varrho \in {\rm Sp}^{k}(M)$ and any $\ga \in \Ga_{\x}$, by
the action of $\ga$ on $\ph_{\psi}$ we have :
$$
(\ga.\ph_{\psi})(\varrho)
=\ph_{\psi}(\ga^{-1}.\varrho)
=\psi(\ga^{-1}.\varrho)(1_{\infty}),
$$
than, since $\psi$ is a homomorphism of $M[G(K_{\infty})]$-modules,
and by the action of $G(K_{\infty})$ on the function
$\psi(\varrho)\in {\frak Aut}_{\x}(L)$ which is right invariant
$\Ga_{\x}$, we have :
$$
\psi(\gamma^{-1}.\varrho)(1_{\infty})
=(\ga^{-1}.\psi(\varrho))(1_{\infty})
=\psi(\varrho)(1_{\infty}\ga^{-1}) =\psi(\varrho)(1_{\infty})
=\ph_{\psi}(\varrho),
$$
we conclude, combining these equalities to preceding ones, that we
have $\ga .\ph_{\psi} =\ph_{\psi}$ for any $\ga \in \Ga_{\x}$.
\smallskip\\
The $M$-linear maps $\Psi_{\x}$ and $\Psi_{\x}^{'}$ are isomorphisms
reciprocal to each other. On one hand we have $\Psi_{\x}\circ
\Psi'_{\x}=I$. Indeed, let $\psi \in {\rm
Hom}_{M[G(K_{\infty})]}({\rm Sp}^{k}(M),\,{\frak Aut}_{\x}(L))$. For
any  $\varrho\in {\rm Sp}^{k}(M)$ and any $g_{\infty}\in
G(K_{\infty})$, We have :
\begin{equation}\label{e5}
\psi_{\varphi_{\psi}}(\varrho)(g_{\infty})
=\varphi_{\psi}(g_{\infty}.\varrho)
=\psi(g_{\infty}.\varrho)(1_{\infty}),
\end{equation}
than, from that fact $\psi$ is a homomorphism of
$M[G(K_{\infty})]$-modules and from that $G(K_{\infty})$ acts on the
function $\psi(g_{\infty})\in {\frak Aut}_{\x}(L)$, we deduce :
$$
\psi(g_{\infty}.\varrho)(1_{\infty})
=[g_{\infty}.\psi(\varrho)](1_{\infty})
=\psi(\varrho)(1_{\infty}g_{\infty}) =\psi(\varrho)(g_{\infty}).
$$
From these equalities and that of (\ref{e5}), we deduce that
$\psi_{\ph_{\psi}}(\varrho)(g_{\infty})=\psi(\varrho)(g_{\infty})$
for any  $\varrho\in {\rm Sp}^{k}(M)$ and any $g_{\infty}\in
G(K_{\infty})$. Consequently, we have
$\psi_{\varphi_{\psi}}=\psi$.\smallskip\\
On the other hand $\Psi'_{\x}\circ \Psi_{\x}=I$. Indeed, let $\ph\in
{\rm Hom}_{M}({\rm Sp}^{k}(M),\, L)^{\Gamma_{\x}}$. Put $\psi_{\ph}=
\Psi_{\x}(\ph)$ and $\ph_{\psi_{\ph}}=\Psi'_{\x}(\psi_{\ph})$.  For
any $\varrho \in {\rm Sp}^{k}(M)$, we have :
$$
\ph_{\psi_{\ph}}(\varrho)
=\psi_{\ph}(\varrho)(1_{\infty})
=\ph(1_{\infty}.\varrho)
=\ph(\varrho).
$$
So $\varphi_{\psi_{\varphi}}=\varphi$. \qed

\section{Automorphic forms and harmonic cochains}

\subsection{A diagram to summarize the isomorphisms seen so far.}

In this paragraph we want to summarize all the preceding results in
a commutative  diagram. This helps in particular to see how to
combine them to get the last result stated in this paper, Theorem
\ref{lien-frmes-autom-cocycles-harm}. To complete the diagram we
need the following lemma :

\begin{lem}\label{diag1} We have an isomorphism of $M$-modules :
$$
\Theta: \; \bigoplus_{\x\in X}{\rm Hom}_{M}(M[\hatI^{k}],
L)^{\Ga_{\x}} \xrightarrow{\cong} {\rm Hom}_{M}(M[\widehat{\mathfrak
I}^{k}], {\cal F})^{G(K)}
$$
which to a family $(h_{\x})_{\x\in X}$ associates
$h:=\Theta((h_{\x})_{\x})$ defined as follows. For any $\si \in
\I^k$, the function $h(\si):G(\A_f)\shortrightarrow L$ is given by :
$$
h(\si)(\g_{f})=h_{\x}(\tau^{-1}\si)
$$
for $\g_{f}=\tau\x\k_{f}\in G(\A_{f})$ with $\tau\in G(K)$ and
$\k_{f}\in \K_{f}$.
\end{lem}
\proof $\;$\\
Let us prove that $\Theta$ is  well defined. First, the definition
doesn't depend on the writing of $\g_f=\tau\x\k_f \in G(\A_f)$.
Indeed, let $\tau'\in G(K)$ and $\k'\in \K_f$ be such that
$\g_f=\tau\x\k_f=\tau'\x\k'_f$ as well. From the equality
$\tau\x\k_f=\tau'\x\k'_f$ we deduce
$\tau^{-1}\tau'=\x\k_{f}\k'^{-1}_{f}\x^{-1}$. Consequently, since
$\k_{f}\k'^{-1}_{f}\in \K_{f}$, we have $\ga =\tau^{-1} \tau' \in
\Ga_{\x}=G(K)\cap \x\K_{f}\x^{-1}$. Therefore, since $h_{\x}$ is
invariant under $\Ga_{\x}$, we have :
$$
h_{\x}(\tau^{-1}\si)=h_{\x}(\tau^{-1}\tau'\tau'^{-1}\si)
=h_{\x}(\ga\tau'^{-1}\si)=(\ga^{-1}.h_{\x})(\tau'^{-1}\si)
=h_{\x}(\tau'^{-1}\si).
$$
Let us prove now that $h$ is invariant under $G(K)$. Let $\tau'\in
G(K)$, $\si\in \widehat{\mathfrak I}^{k}$ and
$\g_{f}=\tau\x\k_{f}\in G({\mathbb A}_{f})$. We have :
\begin{equation}\label{e3}
(\tau'.h)(\si)(\g_{f})=[\tau'.h(\tau'^{-1}\si)](\g_{f})=h(\tau'^{-1}\si)(\tau'^{-1}\g_{f}),
\end{equation}
and observing that $\tau'^{-1}\g_{f}=\tau'^{-1}\tau\x\k_{f}$, we
deduce that :
$$
h(\tau'^{-1}\si)(\tau'^{-1}\g_{f})
=h_{\x}((\tau'^{-1}\tau)\tau'^{-1}\si)
=h_{\x}(\tau^{-1}\si)
=h(\si)(\g_{f}).
$$
From these equalities and that of (\ref{e3}), we get $(\tau'
.h)(\si)(\g_{f})=h(\si)(\g_{f})$; consequently $\tau'.h=h$.

The inverse map sends $h\in{\rm Hom}_{M}(M[\hatI^{k}],{\cal
F})^{G(K)}$ to $(h_{\x})_{\x}$ which is  given by
$$
h_{\x}(\si)=h(\si)(\x).
$$
Let $\x\in X$, let us check that $h_{\x}$ is invariant under
$\Ga_{\x}$. Let $\ga\in \Ga_{\x}$,
$\ga=(\underline{\ga}_{f},\ga_{\infty})$ with
$\ga=\underline{\ga}_{f}\in G(\A_{f})$ and $\ga=\ga_{\infty}\in
G(K_{\infty})$. For any $\si\in \widehat{\mathfrak I}^{k}$, we have
:
$$
(\ga.h_{\x})(\si)=h_{\x}(\ga_{\infty}^{-1}.\si)=h(\ga_{\infty}^{-1}.\si)(\x)
=h(\ga_{\infty}^{-1}.\si)(\underline{\ga}_{f}^{-1}.\underline{\ga}_{f}\x)
$$
Then,
$$
h(\ga_{\infty}^{-1}.\si)(\underline{\ga}_{f}^{-1}.\underline{\ga}_{f}\x)
=[\underline{\ga}_{f}.h(\ga_{\infty}^{-1}.\si)](\underline{\ga}_{f}\x)
=(\ga.h)(\si)(\underline{\ga}_{f}\x) =h(\si)(\x) =h_{\x}(\si).
$$
It is easy to see that this map is the inverse map of $\Theta$, and
so $\Theta$ is an isomorphism de $M$-modules.

\begin{rem}\label{diag2} Let $(h_{\x})_{\x}$ and $h$ be as in the lemma above.
It is easy to check that for any $i$,
$1\leq i\leq 4$, the property {\bf (HCi)} in the definition of
harmonic cochains, cf. Defition \ref{cocycles-harmoniques}, is
satisfied by $h$ if and only if it is by $h_{\x}$ for any $\x$.
\end{rem}

Now we can say that we have the following commutative diagram :

$$
\begin{array}{ccc}
\displaystyle\bigoplus_{\x \in X}{{\frak Harm}}^{k}(M,L)^{\Ga_{\x}}
&
\xrightarrow{\Theta} & {{\frak Harm}}^{k}(M,{\cal F})^{G(K)} \\
{} & {} & {} \\
(\Phi_{\x})_{\x} \;\; \downarrow & {} & \downarrow \, \Phi \\
{} & {} & {} \\
\displaystyle\bigoplus_{\x\in X} {\rm Hom}_{M}({\rm
Sp}^{k}(M),\,L)^{\Ga_{\x}}
& {} & {\rm Hom}_{M}({\rm Sp}^{k}(M),\, {\cal F})^{G(K)}  \\
{} & {} & {} \\
(\Psi_{\x})_{\x} \;\; \downarrow & {} & \downarrow \, \Psi \\
{} & {} & {} \\
\displaystyle\bigoplus_{\x\in X} {\rm Hom}_{M[G(K_{\infty})]}({\rm
Sp}^{k}(M),\, {\mathfrak Aut}_{\x}(L)) & \xrightarrow{\Xi} & {\rm
Hom}_{M[G(K_{\infty})]}({\rm Sp}^{k}(M),\, {\mathfrak
Aut}^{\K_f}(L))
\end{array}
$$
The arrows are isomorphisms of $M$-modules. Indeed, $\Theta$ is a
isomorphism given by Lemma \ref{diag1} and Remarque \ref{diag2}
above, $\Xi$ is an isomorphism by Proposition \ref{frm.aut.infty},
The  $(\Psi_{\x})_{\x}$'s and $\Psi$ are isomorphisms by Proposition
\ref{frm.aut.repr.spec.} and finally the $(\Phi_{\x})_{\x}$'s and
$\Phi$ are isomorphisms by Theorem \ref{maintheorem}.

\subsection{Automorphic forms and harmonic cocycles.}

For any integer $k$, $0\leq k\leq n$, for any $\x\in X$, we denote
by ${{\frak Harm}_{!}^{k}(M,L)}^{\Ga_{\x}}$ the set of harmonic
cochains of degree $k$ and with finite supports modulo $\Ga_{\x}$.
\begin{theo}\label{lien-frmes-autom-cocycles-harm}
Let $M$ be a commutative ring and $L$ be an integral $M$-algebra of
characteristic zero; for every $k$, $0\leq k \leq n$, we have
$M$-isomorphisms :
$$
\bigoplus_{\x \in X} {{\frak Harm}^{k}(M,L)}^{\Ga_{\x}} \cong
{\mbox{Hom}}_{M[G(K_{\infty})]}({\mbox{Sp}}^{k}(M),\; {\frak
Aut}^{\K_f}(L))
$$
and
$$
\bigoplus_{\x \in X} {{\frak Harm}_{!}^{k}(M,L)}^{\Ga_{\x}} \cong
{\mbox{Hom}}_{M[G(K_{\infty})]}({\mbox{Sp}}^{k}(M),\; {\frak
Aut}_{\circ}^{\K_f}(L)).
$$
These isomorphisms are functorial on $\K_f$ (i.e. compatible with
the inclusions  $\K'_f \subseteq \K_f$).
\end{theo}
\proof $\;$\\
The first isomorphism is already given by the preceding diagram,
that is $\Psi \circ \Phi \circ \Theta$ or equally $\Xi \circ
(\Psi_{\x})_{\x}\circ (\Phi_{\x})_{\x}$. We use this last
formulation to prove the second isomorphism. Indeed, by Proposition
\ref{frm.aut.infty}, $\Xi$ induces an isomorphism :
$$
\bigoplus_{\x \in X} {\rm Hom}_{M[G(K_{\infty})]}({\rm Sp}^{k}(M),\,
{\mathfrak Aut}_{\x,\circ}(L)) \cong  {{\rm
Hom}}_{M[G(K_{\infty})]}({\rm Sp}^{k}(M),\; {\frak
Aut}_{\circ}^{\K_{f}}(L)),
$$
so we need only to prove that, for each $\x\in X$,  by the
isomorphism
$$
\Psi_{\x}\circ \Phi_{\x} : {\frak Harm}^{k}(M,L)^{\Ga_{\x}}
\rightarrow {\rm Hom}_{M[G(K_{\infty})]}({\rm Sp}^{k}(M),\,
{\mathfrak Aut}_{\x}(L))
$$
we have $\h\in {\frak Harm}_{!}^{k}(M,L)^{\Ga_{\x}}$ if and only if
its image $\psi=\Psi_{\x}\circ \Phi_{\x}(\h)$ verifies that
$\psi(\varrho) \in {\mathfrak Aut}_{\x,\circ}(L)$ for any $\varrho
\in {\rm Sp}^{k}(M)$. This last assertion is equivalent to
$\psi(\chi_{B_{J}P_{J}})\in {\frak Aut}_{\x,\circ}(L)$; indeed,
$\psi$ is a homomorphism of $M[G(K_{\infty})]$-modules and ${\rm
Sp}^{k}(M)$ as such is generated by $\varrho=\chi_{B_{J}P_{J}}$.
First, let us prove the "only if" part. Let $\h\in {\frak
 Harm}_{!}^{k}(M,L)^{\Ga_{\x}}$, so $\h$ has a finite support modulo
$\Ga_{\x}$. Hence, there exist finitely many $k$-cells
$\si_{1},\ldots ,\si_{r}\in \widehat{\mathfrak I}^{k}$ such that
${\rm supp}({\h})\subseteq \Ga_{\x}\si_{1}\cup \ldots \cup
\Ga_{\x}\si_{r}$. For every $i$,
 $i=1,\ldots ,r$, set  :
$$
{\mathfrak U}_{i}=\{g\in G(K_{\infty}),\, g\si_{J}=\si_{i}\}.
$$
For any $i$ such that ${\mathfrak U}_{i} \neq \emptyset$, for any
$g_{0}\in {\mathfrak U}_{i}$, we have ${\mathfrak
U}_{i}=g_{0}B_{J}$. Hence, see (\ref{BI}) in
\S\,\ref{parahoric-subgroups}, the ${\mathfrak U}_{i}$'s are compact
modulo $Z_{G}(K_{\infty})$ and then so is their union ${\mathfrak
U}=\cup_{i=1}^{r}{\mathfrak U}_{i}$. Furthermore, by definition, for
any $g_{\infty}\in G(K_{\infty})$ we have :
$$
\psi(\chi_{B_{J}P_{J}})(g_{\infty})={\mathfrak
h}(g_{\infty}\si_{J}),
$$
therefore, we have ${\rm supp}\,(\psi(\chi_{B_{J}P_{J}}))\subseteq \Ga_{\x}{\mathfrak U}$. \\
Now, the "if" part. Take $\psi \in {\rm Hom}_{M[G(K_{\infty})]}({\rm
Sp}^{k}(M),\, {\mathfrak Aut}_{\x,0}(L))$, thus for any $\varrho\in
{\rm Sp}^{k}(M)$, its image $\psi(\varrho)$ has a finite support.
Let us prove that the harmonic cochain $\h$ that corresponds to
$\psi$ by the isomorphism $\Psi_{\x}\circ \Phi_{\x}$ is finitely
supported modulo ${\Ga_{\x}}$. From the proof of Proposition
\ref{frm.aut.infty}, there exists ${\mathfrak U}\subseteq
G(K_{\infty})$ so that the quotient $\Ga_{\x}\backslash {\mathfrak
U}/B$ is finite, and so that ${\rm
supp}[\psi(\chi_{C_{I}P_{J}})]\subseteq {\mathfrak U}$, for any
$I\subseteq \Delta$ such that $|\Delta -I|=k$ (for the definition of
the set $C_{I}$, see (\ref{CI}) in \S\,\ref{defCI}). Take
$g_{1},\ldots ,g_{r}$ in $G(K_{\infty})$ so that we have :
$$
\Ga_{\x}\backslash {\mathfrak U}/B =\{\Ga_{\x}g_{i} B,\; 1\leq i\leq
r\}.
$$
If $\h$ is not finitely supported modulo $\Ga_{\x}$, there would be
 $g\in G(K_{\infty})$ and $I\subseteq \Delta$ such that :

\par (1) $\qquad g\si_{I} \notin \bigcup_{I\subseteq \Delta}\bigcup_{i=1}^{r}\Ga_{\x}g_{i}\si_{I}$\\
and
\par (2) $\qquad {\mathfrak h}(g\si_{I})\neq 0$.  \\
The second assertion above is equivalent to
$\psi(\chi_{C_{I}P_{J}})(g)\neq 0$, thus $g\in {\mathfrak U}$.
Contradiction to the first assertion. \qed

 \bigskip
Y. A\"{\i}t Amrane, Facult\'e de Mathématiques, USTHB, El Alia BP , 16111 Algiers. ALGERIA.\\
e-mail adresse : yacinait@hotmail.com

\end{document}